# A New Super-Twisting Algorithm-Based Sliding Mode Observer Design for Fault Estimation in a Class of Nonlinear Fractional Order Systems

Seyed Mohammad Moein Mousavi, Amin Ramezani, HamidReza Momeni

*Department of Electrical and Computer Engineering, Tarbiat Modares University, Tehran, Iran.*

**Abstract**:

This paper is concerned with fault estimation in a class of nonlinear fractional order systems using a new super twisting algorithm based second order step by step sliding mode observer. Since the existing sliding mode observers are troubled with the chattering phenomenon, here a new observer structure is proposed and finite time convergence of error dynamics is proved using fractional order super twisting algorithm (FSTA). Two numerical examples of chaotic fractional order systems and a comparison with respect to a similar observer justify the effectiveness of the proposed observer

## 1. Introduction

Fault diagnosis has been always an important subject in the industry [1]. Sensor and actuator faults can cause failure and damage in physical systems if they are not detected in appropriate time. Fault detection methods are divided into two main types: data based and model based methods. Indeed, model based fault detection methods have been extended during last three decades. In these methods, sensor and actuator faults are detected through the relations between accessible signals. the most popular model based methods are: parameter estimation, observer design and parity space. [2] provides a survey on observer based fault detection. Also observer-based fault detection for nonlinear systems is discussed in [3]. After awareness of existence of a fault, it is important to estimate its shape and domain as an unknown input.

In presence of unknown input, there are two kind of observers. In the first kind [4-5], unknown input must be rejected and it must not affect the estimation of state variables, so it is considered as a disturbance. But in the second kind besides the estimation of states, even the unknown input is reconstructed. These kind of observers are applied to the diagnosis and fault detection problems where the unknown input is considered as the fault signal.

Fractional calculus as a generalization of classic calculus is a mathematic tool which is recently being used in control engineering [6]. Among recent years it has been proven that some systems can be modeled more accurate using fractional order models, compare to integer order ones [7]. Also we can cite Arneodo's system, the Genesio–Tesi's system [8] and Rössler's system [9] as examples of many integer-order chaotic models which are extended to fractional order models. Nevertheless, in working with fractional order systems and observers, stability analysis is usually a more difficult topic compare to integer order models, which is because of special features of fractional calculus. for example, the chain rule is not valid for fractional derivatives and it causes problems using direct lyapunov stability method. In [10] a method is proposed to overcome this problem. Since some physical systems are better modeled using fractional order models, fault estimation in fractional order models is of major significance.

Among different types of observers, sliding mode observers have found a great attention due to their robustness to noise and uncertainties. Sliding mode observers for fault detection and isolation is addressed in [11]. Synchronization of fractional-order systems using adaptive sliding mode techniques is considered in [12-15]. Observability of the states in nonlinear fractional order systems is discussed in [16] and using a first order sliding mode observer, fault is estimated as an unknown input in such system. But as the observer is a first order one, the chattering problem exists. In [17] a second order sliding mode observer is used for fault estimation in a linear fractional order system. a simple method to design a functional observer for linear fractional-order systems with unknown inputs is presented in [18]. A high-order sliding mode observer is proposed in [19] for the pseudo-state and the unknown input estimation of fractional commensurate linear systems

These works either consider linear fractional order systems or first order sliding mode observers. Hence, to decline chattering effect in state and unknown input estimation in a nonlinear fractional order system, it is decent to use a second order sliding mode observer. In integer order case, [20-21] design such observers based on super twisting algorithm. Based on a method suggested in [22] to solve fractional order differential equations, [23] generalizes the super twisting algorithm to fractional order systems and designs a second order sliding mode for observable form of nonlinear fractional order systems. But in this observer, estimation of unknown input is still of the first order and may not reduce the chattering effect, which is the main cause of using a second order sliding mode observer. The aim of tis paper is to introduce a new structure for second order sliding mode observer to reduce the chattering effect on state and unknown input estimations. For this aim, fault is considered as a separate state in dynamic of the observer.

The rest of this paper is organized as follows: some preliminaries on fractional calculus is introduced in section 2. In section 3, the proposed observer and its stability analyze is presented. in section 4 the effectiveness of the proposed observer is discussed by two numerical examples and simulations. Finally, the conclusion remarks are given in section 5.

## 2. Preliminaries

Let $C[a\ b]$ be the space of continuous functions $f(t)$ on $[a\ b]$ and we mean by $C^k$ the space of real-valued functions $f(t)$ with continuous derivatives up to order $k-1$ such that $f^{(k-1)}(t) \in C[a\ b]$ and $f^i(t)$ is the i-th derivative of $f(t)$.

According to [6] there are three main definitions of fractional order derivatives. The Riemann–Liouville integral operator is described as follows.

**Definition 1:** The Riemann–Liouville fractional derivative of order $\alpha$ of $f(t) \in C^m[a\ b]; t \in [a\ b]$:

$$^{RL}_{a}D^{\alpha}_{t}f(t) = \frac{1}{\Gamma(m-\alpha)}\frac{d^m}{dt^m}\int_{a}^{t}(t-\tau)^{m-\alpha-1}f(\tau)d\tau \qquad (1)$$

Where and $\Gamma(\ )$ is Euler's gamma function and it is defined as:

$$\Gamma(\alpha) = \int_{0}^{\infty} v^{\alpha-1}e^{-v}dv \qquad (2)$$

**Definition 2:** Caputo's derivative of order $\alpha$ of $f(t) \in C^m[a\ b]; t \in [a\ b]$:

$$^{C}_{a}D^{\alpha}_{t}f(t) = \frac{1}{\Gamma(m-\alpha)}\int_{a}^{t}(t-\tau)^{m-\alpha-1}\frac{d^m f(\tau)}{dt^m}d\tau \qquad (3)$$

**Definition 3:** Grunwald-Letnikov definition:

$$^{GL}_{a}D^{\alpha}_{t}f(t) = \lim_{h \to \infty}\frac{1}{\Gamma(a)h^{\alpha}}\sum_{j=0}^{\left[\frac{t-a}{h}\right]}\frac{\Gamma(a+j)}{\Gamma(j+1)}f(t-jh) \qquad (4)$$

Where $\left[\frac{t-a}{h}\right]$ denotes the integer part of $\frac{t-a}{h}$.

The drawback of the first definition is that the initial conditions are in terms of the variable's fractional order derivatives. However, Caputo's definition of fractional derivative needs the initial conditions of the main function and not its fractional derivatives. Therefore, in engineering usages Caputo's definition is more common. for the simplicity in the rest of this paper this notation is used:

$$^{C}_{0}D^{\alpha}_{t}f(t) \triangleq D^{\alpha}f(t) \qquad (5)$$

## 3. Observer Design

Consider this observable form of nonlinear fractional order systems in presence of unknown input [16]:

$$\begin{cases} D^\alpha x_1 = x_2, \dots, D^\alpha x_{n-1} = x_n \\ D^\alpha x_n = a(x) + b(x)f \\ D^\alpha f = W \end{cases} \quad (6)$$

Where $0 < \alpha < 1$, $x \in^{n \times 1}$ is the state vector, $a(x)$ and $b(x)$ are Lipchitz functions.

**Assumption 1:** all the states and their first order derivatives are bounded.

**Assumption 2:** W and $D^\alpha W$ are bounded.

Now assume that $b(x) = 1$, then the system (6) finds this form:

$$\begin{cases} D^\alpha x_1 = x_2, \dots, D^\alpha x_{n-1} = x_n \\ D^\alpha x_n = a(x) + f(t) \end{cases} \quad (7)$$

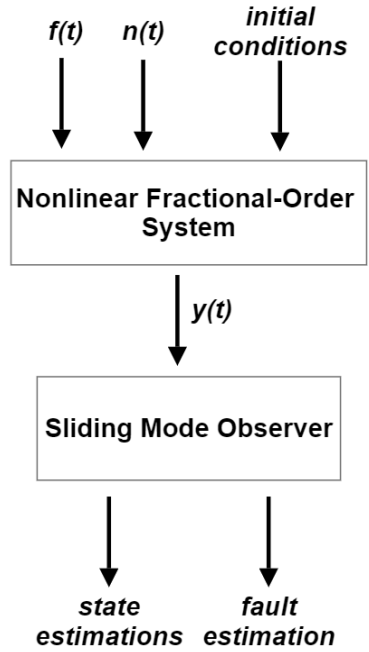

**Figure 1.** general block diagram of the system, observer and their inputs and outputs

Some fractional order chaotic systems such as Arneodo's system and Genesio-Tesi's System are described in form (7). Figure 1 represents a general block diagram of the system, observer and their inputs and outputs. Where $f(t)$ is the fault signal to be estimated, $n(t)$ is noise and $y(t) = x_1(t)$ is system output. A second order observer is proposed in [23] and it is described by these equations:

$$\begin{cases} D^\alpha \hat{x}_1 = \tilde{x}_2 + \lambda_1 |e_1|^{0.5} sign(e_1) \\ D^\alpha \tilde{x}_2 = \alpha_1 sign(e_1) \\ D^\alpha \hat{x}_2 = E_1[\tilde{x}_3 + \lambda_2 |e_2|^{0.5} sign(e_2)] \\ D^\alpha \tilde{x}_3 = E_1 \alpha_2 sign(e_2) \\ D^\alpha \hat{x}_3 = E_2[\tilde{x}_4 + \lambda_3 |e_3|^{0.5} sign(e_3)] \\ \vdots \\ D^\alpha \tilde{x}_{n-1} = E_{n-3} \alpha_{n-2} sign(e_{n-2}) \\ D^\alpha \hat{x}_{n-1} = E_{n-2}[\tilde{x}_n + \lambda_{n-1} |e_{n-1}|^{0.5} sign(e_{n-1})] \\ D^\alpha \tilde{x}_n = E_{n-2} \alpha_{n-1} sign(e_{n-1}) \\ D^\alpha \hat{x}_n = E_{n-1}[\tilde{\theta} + \lambda_n |e_n|^{0.5} sign(e_n)] \\ D^\alpha \tilde{\theta} = E_{n-1} \alpha_n sign(e_n) \end{cases} \quad (8)$$

This structure does not consider the fault as a state in observer dynamics. Here we propose this structure for second order observer:

$$\begin{cases} D^\alpha \hat{x}_1 = \tilde{x}_2 + \lambda_1 |e_1|^{0.5} sign(e_1) \\ D^\alpha \tilde{x}_2 = \alpha_1 sign(e_1) \\ D^\alpha \hat{x}_2 = E_1[\tilde{x}_3 + \lambda_2 |e_2|^{0.5} sign(e_2)] \\ D^\alpha \tilde{x}_3 = E_1 \alpha_2 sign(e_2) \\ D^\alpha \hat{x}_3 = E_2[\tilde{x}_4 + \lambda_3 |e_3|^{0.5} sign(e_3)] \\ \vdots \\ D^\alpha \tilde{x}_{n-1} = E_{n-3} \alpha_{n-2} sign(e_{n-2}) \\ D^\alpha \hat{x}_{n-1} = E_{n-2}[\tilde{x}_n + \lambda_{n-1} |e_{n-1}|^{0.5} sign(e_{n-1})] \\ D^\alpha \tilde{x}_n = E_{n-2} \alpha_{n-1} sign(e_{n-1}) \\ D^\alpha \hat{x}_n = E_{n-1}[a(\tilde{x}) + \tilde{f} + \lambda_n |e_n|^{0.5} sign(e_n)] \\ D^\alpha \tilde{f} = E_{n-1} \alpha_n sign(e_n) \\ D^\alpha \hat{f} = E_n [\tilde{\theta} + \lambda_{n+1} |e_f|^{0.5} sign(e_f)] \\ D^\alpha \tilde{\theta} = E_n \alpha_{n+1} sign(e_f) \end{cases} \quad (9)$$

And we define Estimation errors as:

$$e_i = \tilde{x}_i - \hat{x}_i \quad (10)$$

With $\tilde{x}_1 = x_1$ and $x_1$ is the output which is the only accessible signal. We have:

$$E_i = 1 \; if \; |e_j| = |\tilde{x}_j - \hat{x}_j| \leq \varepsilon \; for \; all \; j \leq i \; else \; E_i = 0 \quad (11)$$

$\varepsilon$ is a small positive scalar and the observer gains are all positive coefficients.

**Theorem 1**: (fractional order super twisting algorithm(FSTA)) [23]

The structure of FSTA is described by the following equations:

$$\begin{cases} D^\alpha \xi_1 = \xi_2 - \lambda_1 |\xi_1|^{0.5} sign(\xi_1) \\ D^\alpha \xi_2 = -\alpha_1 sign(\xi_1) + \rho(t) \end{cases} \quad (12)$$

Where $0 < q < 1$ and $\xi_1(t) \in R$ and $\xi_2(t) \in R$ are scalar real variables. $\rho(t)$ represents the perturbation term assumed to be bounded. The coefficients $\lambda_1 > 0$ and $\alpha_1 > 0$ are the tuning strictly positive gains of the FSTA

Suppose that the perturbation term $\rho(t)$ is bounded that is $|\rho(t)| \leq L$. Then, there exist gains $\lambda_1, \alpha_1$ such that for any initial conditions $\xi_1(0)$ and $\xi_2(0)$, the variables $\xi_1(t)$ and $\xi_2(t)$, solutions of the FSTA (12) converge to zero in finite time $T_s$:

$$T_s = (\Gamma(\alpha + 1) v_s)^{\frac{1}{\alpha}} \quad (13)$$

Where $v_s$ is the time convergence of the equivalent integer-order system of (11).

**Theorem 2**: Consider the system (7) and the observer (9). Assume that the system is bounded-input-bounded-state, then for any bounded state $x_0, \hat{x}_0$ and any bounded fault f(t), there exist constants $\lambda_i, \alpha_i$ such that the observer state converges to real states in finite time and the estimated fault converges to the fault in finite time.

**Proof.** The proof comprises of the following steps:

Step1: $E_i = 0 \; for \; all \; i$

The estimation errors are driven by the following equations:

$$\begin{cases} D^\alpha e_1 = x_2 - \tilde{x}_2 - \lambda_1 |e_1|^{0.5} sign(e_1) \\ D^\alpha \tilde{x}_2 = \alpha_1 sign(e_1) \end{cases} \quad (14)$$

These equations are of the form (12) with $\rho_1 = x_3$. So using FSTA we can say that $e_1$ and its fractional derivatives converge to zero in finite time and Considering (14), $x_2$ converges to $\tilde{x}_2$.

Step 2: $E_1 = 1$ and:

$$D^\alpha e_1 = 0$$
$$D^\alpha e_2 = x_3 - \tilde{x}_3 - \lambda_2 |e_2|^{0.5} sign(e_2) \qquad (15)$$
$$D^\alpha \tilde{x}_3 = \alpha_2 sign(e_2)$$

Similarly, a fractional order super twisting algorithm appears with $\rho_2 = x_4$. Hence, $e_2$ and its derivatives would be 0 and $\tilde{x}_3$ converges to $x_3$.

⋮

Step n: $E_1 = E_2 = \cdots = E_{n-1} = 1$, $E_n = 0$, $x_i = \tilde{x}_i$ for all i

$$\begin{aligned}D^\alpha e_1 &= \cdots = D^\alpha e_{n-1} = 0 \\ D^\alpha e_n &= a(x) + f(t) - a(\tilde{x}) - \tilde{f}(t) - \lambda_n |e_n|^{0.5} sign(e_n) \\ D^\alpha \tilde{f} &= \alpha_n sign(e_n)\end{aligned} \qquad (16)$$

Since $\tilde{x}$ is converged to $x$, we have $a(x) = a(\tilde{x})$. So we can write:

$$\begin{aligned}D^\alpha e_n &= f(t) - \tilde{f}(t) - \lambda_n |e_n|^{0.5} sign(e_n) \\ D^\alpha \tilde{f} &= \alpha_n sign(e_n) \\ D^\alpha (f - \tilde{f}) &= W - \alpha_n sign(e_n)\end{aligned} \qquad (17)$$

Which is again structure of FSTA with $\rho_n = W$. So we can say $e_n$ and $\tilde{f}$ converge to 0 and $f$ respectively.

Step n+1: $E_1 = E_2 = \cdots = E_{n-1} = E_n = 1$, $f(t) = \tilde{f}(t)$

$$\begin{aligned}D^\alpha e_1 &= \cdots = D^\alpha e_n = 0 \\ D^\alpha e_f &= W - \tilde{\theta} - \lambda_{n+1} |e_f|^{0.5} sign(e_f) \\ D^\alpha \tilde{\theta} &= \alpha_{n+1} sign(e_f)\end{aligned} \qquad (18)$$

So we can write:

$$\begin{aligned}D^\alpha e_f &= W - \tilde{\theta} - \lambda_{n+1} |e_f|^{0.5} sign(e_f) \\ D^\alpha (W - \tilde{\theta}) &= -\alpha_{n+1} sign(e_f) + D^\alpha W\end{aligned} \qquad (19)$$

Here $\rho_n = D^\alpha W$ and is assumed to be bounded. So $e_f$ converges to 0 and the proof is completed.

**Remark 1**: in systems of the form (6) where $b(x) \neq 1$, one can consider the term $b(x)f(t)$ as an additive fault say $D$. It was proved that observer (9) estimates this kind of unknown input in finite time. So an estimation of $f(t)$ in such system can be:

$$\hat{f}(t) = b(\hat{x})^{-1}\hat{D} \tag{20}$$

**Remark 2**: with use of observer (8) proposed in [23], the unknown input is not considered as a state in the observer structure, and the final estimation of the fault is:

$$\hat{f}(t) = b^{-1}(\tilde{x})(\tilde{\theta} - a(\tilde{x})) \tag{21}$$

This estimation only depends on $\tilde{x}_i$ and $\tilde{\theta}$. Despite all the $\tilde{x}_i$, $\hat{x}_i$ converge to their real states, but $\tilde{x}_i$ have a first order dynamic and it is more affected by chattering. Example 2 explains this remark better.

## 4. Numerical simulations

**Example 1**: In this example we consider Arneodo's system described in [8] and we illustrate the effectiveness of proposed observer through fault estimation in this chaotic nonlinear fractional order system in presence of state noise. Consider:

$$\begin{cases} D^\alpha x_1 = x_2 \\ D^\alpha x_2 = x_3 \\ D^\alpha x_3 = -\beta_1 x_1 - \beta_2 x_2 - \beta_3 x_3 + \beta_4 x_1^3 + f(t) + n(t) \\ f(t) = 0.4\cos(t) \end{cases} \tag{22}$$

Where $\alpha = 0.97$. $f(t)$ is the fault to be estimated and $n(t)$ is a white Gaussian noise with zero mean and variance of 1.5. Here the system parameter are chosen as $\beta_1 = -5.5, \beta_2 = 3.5, \beta_3 = 0.8, \beta_4 = -1.0$. Figure 2 shows chaotic behavior of the system in absence of fault and noise and computational time 200s, for time step $h = 0.001$ and initial conditions $x_{1_0} = -0.2; x_{2_0} = 0.5; x_{3_0} = 0.2$

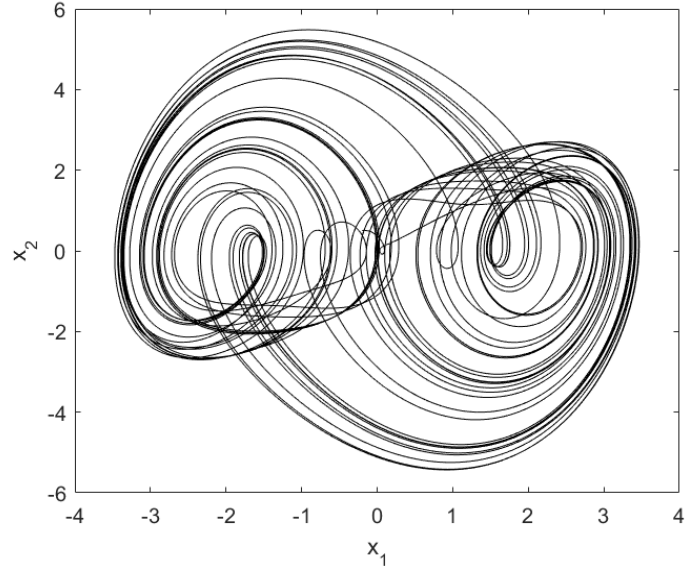

**Figure 2.** chaotic behavior of system in example 1

Now we consider this structure for the observer:

$$\begin{cases} D^\alpha \hat{x}_1 = \tilde{x}_2 + \lambda_1 |e_1|^{0.5} sign(e_1) \\ D^\alpha \tilde{x}_2 = \alpha_1 sign(e_1) \\ D^\alpha \hat{x}_2 = E_1[\tilde{x}_3 + \lambda_2 |e_2|^{0.5} sign(e_2)] \\ D^\alpha \tilde{x}_3 = E_1 \alpha_2 sign(e_2) \\ D^\alpha \hat{x}_3 = E_2[-\beta_1 x_1 - \beta_2 \tilde{x}_2 - \beta_3 \tilde{x}_3 + \beta_4 x_1^3 + \tilde{f} + \lambda_3 |e_3|^{0.5} sign(e_3)] \\ D^\alpha \tilde{f} = E_2 \alpha_3 sign(e_3) \\ D^\alpha \hat{f} = E_3 \left[\tilde{\theta} + \lambda_4 |e_f|^{0.5} sign(e_f)\right] \\ D^\alpha \tilde{\theta} = E_3 \alpha_4 sign(e_f) \end{cases} \quad (23)$$

We choose observer gains as $\lambda_1 = 1$; $\alpha_1 = 10$; $\lambda_2 = 1$; $\alpha_2 = 200$; $\lambda_3 = 10$; $\alpha_3 = 50$; $\lambda_4 = 100$; $\alpha_4 = 100$; $x_1$ as the output of the system is the only accessible signal and it is the observer input. The initial conditions of the observer parameters are zero. It should be noted that observer gains are chosen by trial and error. If these gains are chosen small then the estimation errors will not converge to zero, also if we choose big values for the gains, the chattering problem will appear and the estimation error magnitudes increase. Therefore, a future work can be optimum choice of these gains. Also, it should be mentioned that the smaller time step h is chosen, the less chattering will appear in observer estimations, but the simulation time will increase as well.

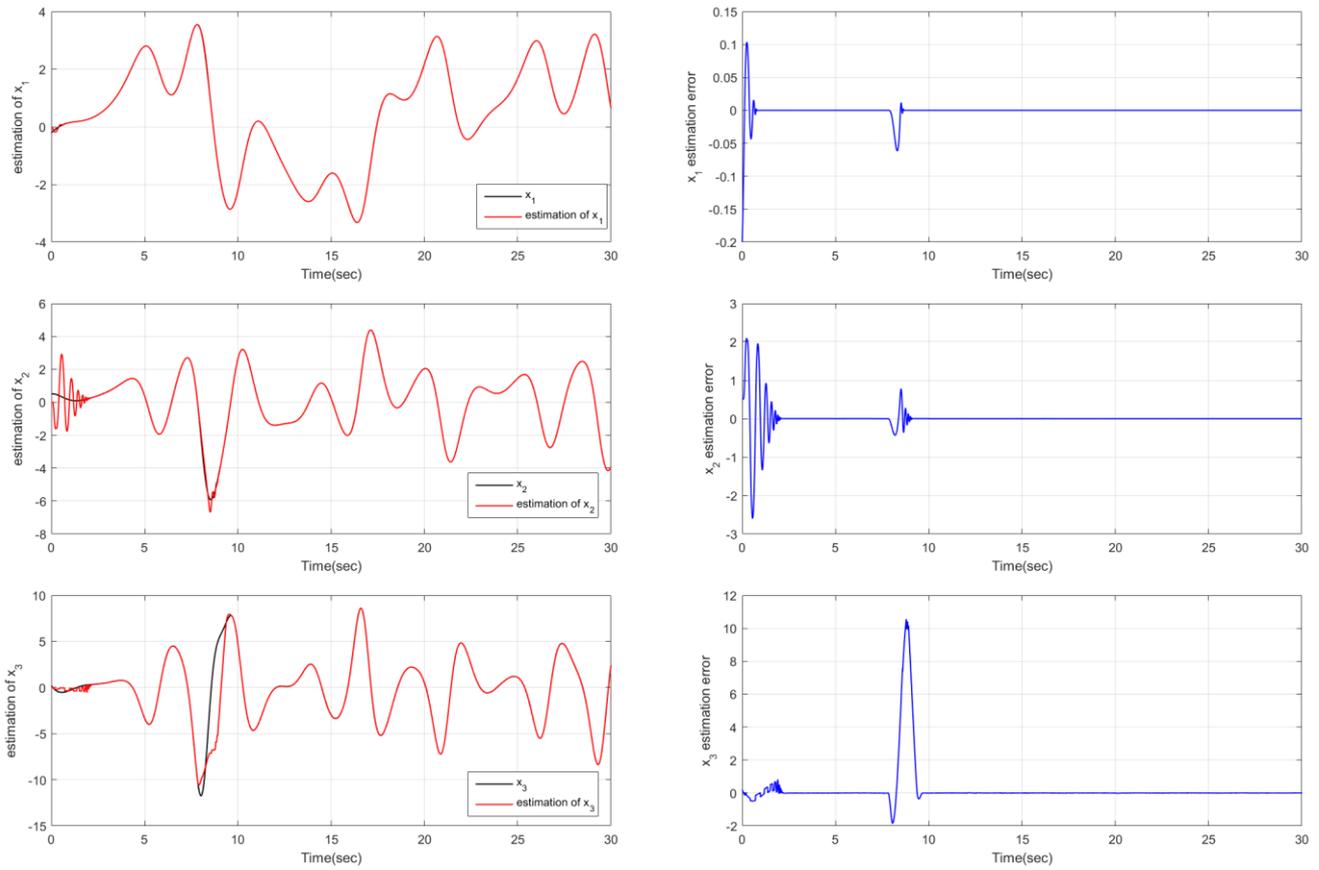

**Figure 3.** state estimations in example 1

Simulations are performed using Grünwald–Letkinov's definition of the fractional derivative. State estimations and their errors are depicted in Figure 3. Also Figures 4 represents the fault signal and its estimation and the estimation error. One can note that the fault signal is reconstructed and the estimation errors converge to zero. Hence this simulation endorses the effectiveness of the proposed observer.

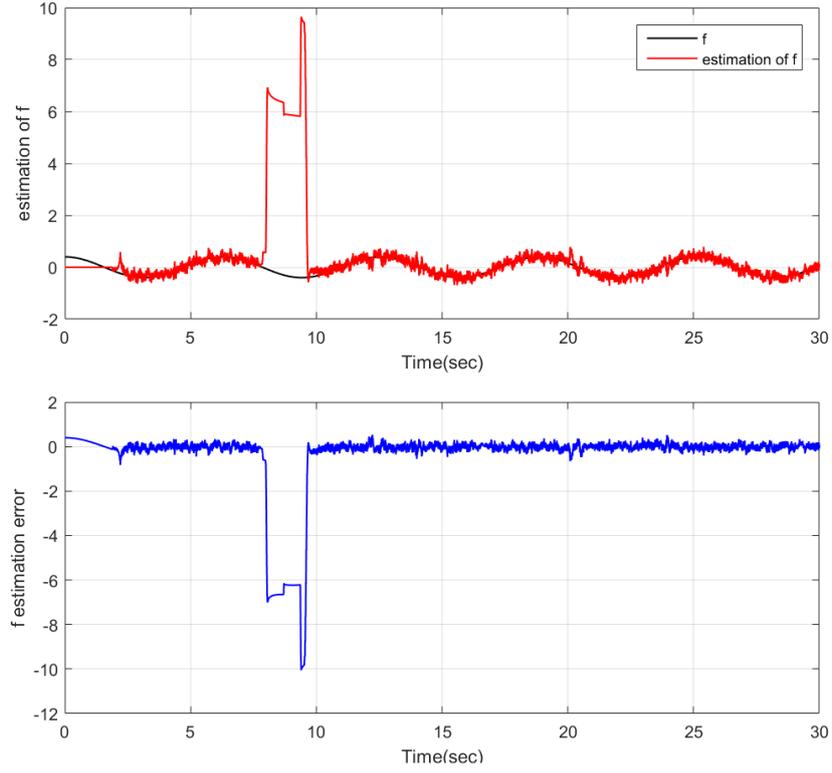

**Figure 4.** fault estimation in example 1

**Example 2:** In this example we consider Genesio-Tesi's System described in [8] and we compare performances of the proposed observer (9) and the similar observer (8) introduced in [23]. Consider:

$$\begin{cases} D^\alpha x_1 = x_2 \\ D^\alpha x_2 = x_3 \\ D^\alpha x_3 = -\beta_1 x_1 - \beta_2 x_2 - \beta_3 x_3 + \beta_4 x_1^2 + f(t) \\ f(t) = 0.06 \sin(t) \end{cases} \quad (24)$$

Where $\alpha = 0.9$. Again f(t) is the fault to be estimated. for the following parameters: $\beta_1 = 1, \beta_2 = 1.1, \beta_3 = 0.44, \beta_4 = 1.0.$

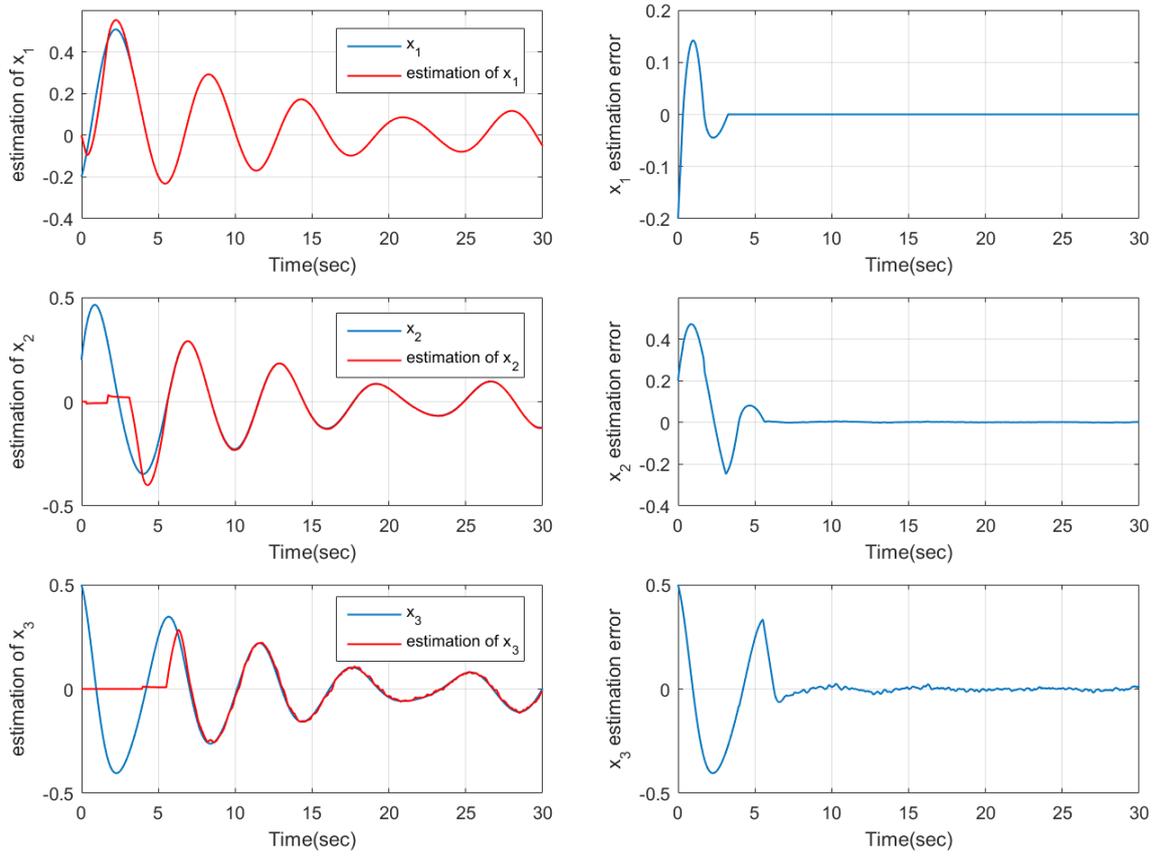

**Figure 5**. state estimations and their errors in example 2 using new observer

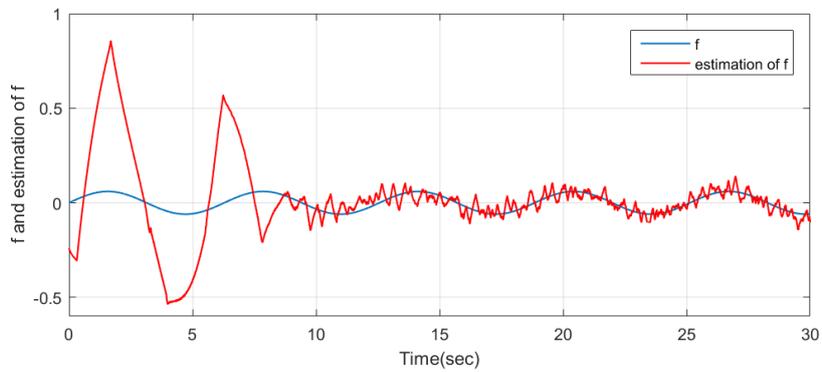

**Figure 6**. fault estimation in example 2 using observer (8)

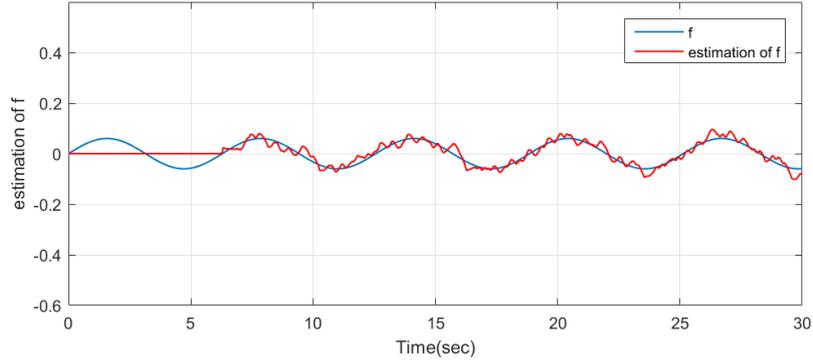

**Figure 7.** fault estimation in example 2 using observer (9)

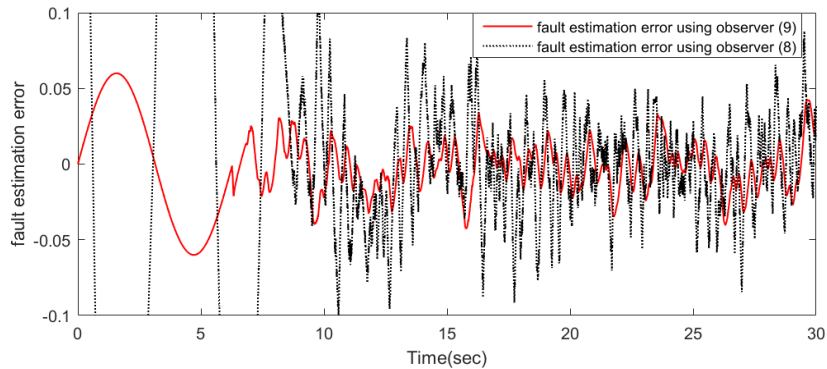

**Figure 8.** fault estimation errors in example 2

As the results are notably sensitive to the gains, here we chose the same gains for both observers. So all the gains are adopted equal to 0.5. In state estimation, both observers show a good performance, Figure 5 shows for example, state estimations using the new observer. But there is a remarkable difference in fault estimation as we expected. Figures 6-8 depict that using the new observer, the fault estimation error magnitude and also its chattering are declined compare to the similar observer.

## 5. Conclusion

In this paper, a new structure was introduced for second order step by step sliding mode observer, in order to estimate states and fault as an unknown input in a class of nonlinear fractional order systems. based on fractional order super twisting algorithm, stability and finite time convergence of error dynamics was proved. Finally, two numerical simulations were performed to illustrate the effectiveness of the proposed observer. A comparison with a similar observer in terms of unknown input estimation depicted that the new observer is less affected by chattering phenomenon. Since the estimations are hardly sensitive to the observer gains, a further study here can be optimum choice of observer gains.